\documentclass[10pt,journal]{IEEEtran}
\usepackage{amsmath,amssymb,epsfig}
\usepackage{longtable}
\allowdisplaybreaks

\begin{document}
\pagenumbering{gobble}
\newcommand*{\AuthorFont}{%
      \fontsize{23}{28}%
      \selectfont}
\title{\AuthorFont Cooperative Compressive Power Spectrum Estimation\vspace{-1ex}} 
\author{Dyonisius Dony Ariananda$^\ast$, Daniel Romero$^\dag$, and Geert Leus$^\ast$\\
$^\ast$Faculty of EEMCS, Delft University of Technology, The Netherlands\\
$^\dag$Dept. of Signal Theory and Communications, University of Vigo, Spain\\
\{d.a.dyonisius, g.j.t.leus\}@tudelft.nl and dromero@gts.uvigo.es.\vspace{-4ex}} 
\maketitle
\thispagestyle{empty}
\begin{abstract}
We examine power spectrum estimation from wide-sense stationary signals received at different wireless sensors. 
We organize multiple sensors into several groups, where each group estimates the temporal correlation only at particular lags, which are different from group to group. A fusion centre collects all the correlation estimates from different groups of sensors, and uses them to estimate the power spectrum. This reduces the required sampling rate per sensor. 
We further investigate the conditions required for the system matrix to have full column rank, 
which allows for a least-squares reconstruction method. 
\end{abstract}
\vspace{-3mm}
\section{Introduction}\label{introduction}

The emergence of compressive sampling has renewed the interest in spectral analysis. The work in~\cite{Eldar}, for instance, focuses on 
signal reconstruction from sub-Nyquist rate samples produced by a multi-coset sampler. However, when only the statistics of the received signal are of interest (such as in a cognitive radio (CR) application), attempting to reconstruct the uncompressed waveform is unnecessary. 
In this case, power spectrum reconstruction from compressive measurements becomes more appropriate~\cite{TSP12}~\cite{XiaodongWang}. 
Compressive power spectrum reconstruction for a wide-sense stationary (WSS) signal is possible, even for a non-sparse power spectrum, 
by exploiting 
the Toeplitz structure of its 
correlation matrix~\cite{TSP12}. Unlike~\cite{TSP12},~\cite{XiaodongWang} 
considers 
a multiband signal (which is not necessarily WSS), where the spectra at different bands are uncorrelated. This allows~\cite{XiaodongWang} to exploit the diagonal structure of the correlation matrix of the entries at different bands. 

In wireless communications, the received user signal might suffer from fading 
and the use of a single 
receiver to compressively reconstruct either the spectrum or the power spectrum of the user signal might be insufficient to reach the required performance. 
In order to 
exploit channel diversity,~\cite{FZheng} 
proposes a 
cooperative compressive wideband spectrum sensing method for CR networks, 
which 
also 
reduces 
the required sampling rate per sensor. 
However, the aim 
to reconstruct the spectrum or the spectrum support 
requires the original spectrum to be sparse. 
This inspired~\cite{AsilomarCR} 
to extend the power spectrum estimation method of~\cite{TSP12} 
into a cooperative scenario. In~\cite{AsilomarCR}, the exploitation of the cross-spectra between the compressive measurements 
at 
different sensors 
reduces the required sampling rate per sensor without requiring the power spectrum to be sparse, 
but it builds upon the knowledge of 
the channel state information (CSI). 
As in~\cite{AsilomarCR}, our work focuses on cooperative compressive power spectrum 
estimation 
but we do not need CSI as we do not exploit the cross-spectra between 
different sensors. 
We consider groups of sensors where different groups employ different sub-Nyquist sampling patterns. Each group estimates the temporal correlation only for certain lags and not for the entire 
correlation support. The fusion centre (FC) collects the correlation estimates 
at different lags produced by different groups of sensors.
The combined 
correlation values at the FC, which should include all the lags in the considered correlation support, are then used to estimate the power spectrum. The required sampling rate per sensor can thus be lower than 
in the single sensor case presented in~\cite{TSP12}, while the channel diversity can still be exploited. 

{\it Notation:} Upper (lower) boldface letters are used to denote matrices (column vectors).
\section{Theoretical Model}\label{system_model}

Consider $Z$ groups of $P$ wireless sensors sensing time-domain WSS signals. At the $(p+1)$-th sensor of the $(z+1)$-th group, we collect $\tilde{N}$ 
Nyquist-rate samples, split them into $L$ blocks of $N={\tilde{N}}/{L}$ consecutive samples, and collect the samples in the $(l+1)$-th block as ${\bf x}_{z,p}[l]=[x_{z,p}[lN],x_{z,p}[lN+1],\dots,x_{z,p}[lN+N-1]]^T$, 
with $l=0,1,\dots,L-1$, $p=0,1,\dots,P-1$, $z=0,1,\dots,Z-1$, and $x_{z,p}[\tilde{n}]$ 
the $(\tilde{n}+1)$-th sample. We collect the $(n+1)$-th indices from each block into the set 
$\{\tilde{n} \in\{0,1,\dots,\tilde{N}-1\}|\tilde{n}\text{ mod }N=n\}$ 
and label this set as the {\it $(n+1)$-th coset}, with the {\it coset index} of the $(n+1)$-th coset given by $n$. 
It is hence possible to view the complete set of $\tilde{N}$ 
samples as the output of a 
multi-coset sampler~\cite{Eldar} with $N$ cosets and $L$ samples per coset.
Next, we introduce compression in each sensor 
by defining one unique sub-Nyquist 
sampling pattern for each group. 
For the $(z+1)$-th group, we define an $M\times N$ selection matrix ${\bf C}_z$, whose rows are selected from the rows of the $N \times N$ identity matrix ${\bf I}_N$.
The $\frac{M}{N}$-rate compressed version of ${\bf x}_{z,p}[l]$ is written as
\vspace{-0.5mm}
\begin{equation}
{\bf y}_{z,p}[l]={\bf C}_z{\bf x}_{z,p}[l],\:\:z=0,\dots,Z-1,\:\:p=0,\dots,P-1,
\label{eq:comp_correlogram}
\vspace{-0.5mm}
\end{equation}
with ${\bf y}_{z,p}[l]=[{y}^{(0)}_{z,p}[l],{y}^{(1)}_{z,p}[l],\dots,{y}^{(M-1)}_{z,p}[l]]^T$. 
We assume that all $\{{\bf C}_z\}_{z=0}^{Z-1}$ in~\eqref{eq:comp_correlogram} 
have $M$ rows to simplify the discussion though it is possible for every ${\bf C}_z$ to have 
a different number of rows. By writing ${\bf C}_z$ in~\eqref{eq:comp_correlogram} in terms of its rows, i.e., 
${\bf C}_z=[{\bf c}^{(0)}_{z},{\bf c}^{(1)}_{z},\dots,{\bf c}^{(M-1)}_{z}]^T$, 
with ${\bf c}^{(m)}_{z}=[{c}^{(m)}_{z}[0],$ ${c}^{(m)}_z[-1],\dots,{c}^{(m)}_{z}[1-N]]^T$, 
we can write 
\vspace{-1.5mm}
\begin{equation}
{y}^{(m)}_{z,p}[l]=
\sum_{n=1-N}^{0}c^{(m)}_{z}[n]{x}_{z,p}[lN-n], 
\label{eq:ympz_filter}
\vspace{-1.5mm}
\end{equation}
for $m=0,1,\dots,M-1$. As~\eqref{eq:ympz_filter} can be viewed as a filtering operation of a WSS sequence ${x}_{z,p}[\tilde{n}]$ by a filter $c^{(m)}_{z}[n]$ followed by an $N$-fold decimation, 
$\{{y}^{(m)}_{z,p}[l]\}_{m=0}^{M-1}$ 
forms a set of jointly WSS sequences. 
We then define the set of indices of the $M$ selected cosets in~\eqref{eq:comp_correlogram} (corresponding to the set of indices of the $M$ rows of ${\bf I}_N$ used to form ${\bf C}_z$) as 
\vspace{-0.5mm}
\begin{equation}
\mathcal{M}_z=\{n^{(0)}_{z}, n^{(1)}_{z}, \dots, n^{(M-1)}_{z}\},
\label{eq:mathcalMz}
\vspace{-1mm}
\end{equation}
with $0\leq n^{(0)}_{z}<n^{(1)}_{z}<\dots< n^{(M-1)}_{z}\leq N-1$. 

We focus on the case where the different sensors observe the same statistics of the received signals, which can be motivated by the fact that they are observing the same user signals, that they experience different random fading, and that they use automatic gain control to compensate the difference in the fading variance (path loss and shadowing) in each user band. 
In this case, we define 
${r}_x[\tilde{n}]=E[{x}_{z,p}[\tilde{n}']{x}^*_{z,p}[\tilde{n}'-\tilde{n}]]$, i.e., $r_x[\tilde{n}]$ does not vary with sensor indices $p$ or group indices $z$. 
Note that~\eqref{eq:mathcalMz} allows us to write $c^{(m)}_{z}[n]$ in~\eqref{eq:ympz_filter} as $c^{(m)}_{z}[n]=\delta[n+n^{(m)}_{z}]$. As a result, the deterministic cross-correlation between $c^{(m)}_{z}[n]$ and $c^{(m')}_{z}[n]$, defined as $r^{(m,m')}_{c_z}[n]=\sum_{n'=1-N}^{0}c^{(m)}_{z}[n']c^{(m')*}_{z}[n'-n]$, is given by
\vspace{-1mm}
\begin{equation}
r^{(m,m')}_{c_z}[n]=\delta[n+n^{(m)}_{z}-n^{(m')}_{z}]
\vspace{-1.5mm}
\label{eq:rczmm}
\end{equation}
and the correlations between the measurements at the different cosets $m$ and $m'$, i.e., 
${r}^{(m,m')}_{{y}_{z}}[l]
=E[{y}^{(m)}_{z,p}[l']{{y}^{(m')*}_{z,p}[l'-l]}]$, 
are related to $r_x[\tilde{n}]$ as
\vspace{-2.5mm}
\begin{align}
{r}^{(m,m')}_{{y}_{z}}[l]&=\sum_{n=1-N}^{N-1}r^{(m,m')}_{c_z}[n]r_x[lN-n] \nonumber \\
\vspace{-1mm}
&=r_x[lN+n^{(m)}_{z}-n^{(m')}_{z}].
\vspace{-2.5mm}
\label{eq:rypzmm'filter}
\end{align}
\hspace*{1mm}
{\it Remark 1: 
This paper focuses on estimating ${r}_x[\tilde{n}]$ in~\eqref{eq:rypzmm'filter} at lags $1-{N}\leq \tilde{n} \leq {N}-1$ from the sample estimates of ${r}^{(m,m')}_{{y}_{z}}[l]$ in~\eqref{eq:rypzmm'filter} with $L\geq 2$. For this purpose, we can find from~\eqref{eq:mathcalMz}-\eqref{eq:rypzmm'filter} that 
we need to consider only ${r}^{(m,m')}_{{y}_{z}}[0]$ for all $m$ and $m'$, ${r}^{(m,m')}_{{y}_{z}}[1]$ for all $m<m'$, and ${r}^{(m,m')}_{y_z}[-1]$ for all $m>m'$.}
\newline
\hspace*{1mm}Considering~\eqref{eq:rczmm}, 
we write $\{{r}^{(m,m')}_{{y}_z}[l]\}_{l=-1}^{1}$ in~\eqref{eq:rypzmm'filter} as 
\begin{subequations}\label{eq:rypzlagall}
\setlength{\abovedisplayskip}{1pt}
\setlength{\belowdisplayskip}{1.5pt}
\begin{align}
&{r}^{(m,m')}_{{y}_z}[0]
=r^{(m,m')}_{c_z}[0]{r}_x[0]+{\bf r}^{(m,m')T}_{c_z}[1]{\bf r}_x[-1]\nonumber \\
&+{\bf r}^{(m,m')T}_{c_{z}}[-1]{\bf r}_x[1],\:\: m,m'=0,1,\dots,M-1,\label{eq:rypzlag0}\\
&{r}^{(m,m')}_{{y}_{z}}[1]={\bf r}^{(m,m')T}_{c_z}[1]{\bf r}_x[1],\nonumber \\
&\quad\quad\quad\quad\quad m=0,1,\dots,M-2,\quad m'>m,\label{eq:rypzlag1}\\
&{r}^{(m,m')}_{{y}_{z}}[-1]={\bf r}^{(m,m')T}_{c_z}[-1]{\bf r}_x[-1], \nonumber \\
&\quad\quad\quad\quad\quad m'=0,1,\dots,M-2,\quad m>m'\label{eq:rypzlagmin1},
\end{align}
\end{subequations}
where we have 
\begin{subequations}\label{eq:rczmlag}
\setlength{\abovedisplayskip}{2pt}
\setlength{\belowdisplayskip}{3pt}
\begin{align}
\small
\label{eq:rxzpmin1}{\bf r}_{{x}}[-1]&=[{r}_{{x}}[1-N],\dots,r_{{x}}[-2],r_{{x}}[-1]]^T, \\ 
\label{eq:rxzp1}{\bf r}_{{x}}[1]&=[{r}_{{x}}[1],r_{{x}}[2],\dots,r_{{x}}[N-1]]^T, \\
\label{eq:rczmlagmin1}{\bf r}^{(m,m')}_{c_z}[-1]&=[r^{(m,m')}_{c_z}[-1],\dots,r^{(m,m')}_{c_z}[1-N]]^T,\\ 
\label{eq:rczmlag1}{\bf r}^{(m,m')}_{c_z}[1]&=[r^{(m,m')}_{c_z}[N-1],\dots,r^{(m,m')}_{c_z}[1]]^T.
\end{align}
\end{subequations}
Observe from~\eqref{eq:rczmm},~\eqref{eq:rczmlagmin1}, and~\eqref{eq:rczmlag1} that the first, the second, and the third terms 
in~\eqref{eq:rypzlag0} are non-zero only if $m=m'$, $m<m'$, and $m>m'$, respectively. 
Based on this fact,~\eqref{eq:rypzlag1},~\eqref{eq:rypzlagmin1}, and the Hermitian property of ${r}^{(m,m')}_{{y}_{z}}[l]$, i.e., ${r}^{(m,m')}_{{y}_{z}}[l]={r}^{(m',m)*}_{{y}_{z}}[-l]$, 
we just need to consider the correlations 
${r}^{(m,m')}_{{y}_{z}}[0]$ for $m\geq m'$ and ${r}^{(m,m')}_{{y}_{z}}[1]$ for $m'>m$, since they contain all relevant information. We then define 
\begin{subequations}\label{eq:Rclags}
\setlength{\abovedisplayskip}{3pt}
\setlength{\belowdisplayskip}{3pt}
\begin{align}  
\label{eq:Rclags0}&{\bf r}_{c_z}[0]=[\dots, r^{(m,m)}_{c_z}[0],\dots]^T,\: m=0,\dots,M-1, \\
\label{eq:Rczlagsmin1}&{\bf R}_{c_z}[-1]=[\dots, {\bf r}^{(m,m')}_{c_z}[-1], \dots]^T, \nonumber \\
&\quad\quad\quad\quad m'=0,1,\dots,M-2,\quad m>m', \\
\label{eq:Rczlags1}&{\bf R}_{c_z}[1]=[\dots, {\bf r}^{(m,m')}_{c_z}[1], \dots]^T,\nonumber \\
&\quad\quad\quad\quad m=0,1,\dots,M-2,\quad m'>m,
\end{align}
\end{subequations}
with ${\bf r}_{c_z}[0]={\bf 1}_M$ (as is clear from~\eqref{eq:rczmm}) an $M\times 1$ vector having ones in all entries and both ${\bf R}_{c_z}[1]$ as well as ${\bf R}_{c_z}[-1]$ $\frac{M(M-1)}{2}\times (N-1)$ matrices. We can then use~\eqref{eq:Rclags} to write 
\begin{subequations}\label{eq:rypzall} 
\setlength{\abovedisplayskip}{3pt}
\setlength{\belowdisplayskip}{3pt}
\begin{align}  
{\bf r}^{(0)}_{{y}_z}[0]&={\bf r}_{c_z}[0]{r}_{{x}}[0]={\bf 1}_M{r}_{{x}}[0],\\
\label{eq:rxpz_one}
\begin{bmatrix} {\bf r}^{(+)}_{{y}_{z}}[0] \\
{\bf r}^{(-)}_{{y}_{z}}[1]
\end{bmatrix}&=
\begin{bmatrix} {\bf R}_{c_z}[-1] \\
{\bf R}_{c_z}[1]
\end{bmatrix}{\bf r}_{{x}}[1]={\bf R}_{c_z}{\bf r}_{{x}}[1],
\end{align}
\end{subequations}
with ${\bf r}^{(0)}_{{y}_z}[0]=[\dots,{r}^{(m,m)}_{{y}_z}[0],\dots]^T$ for $m=0,1,\dots,M-1$, 
${\bf r}^{(+)}_{{y}_z}[0]=[\dots,{r}^{(m,m')}_{{y}_z}[0],\dots]^T$ 
for $m'=0,1,\dots,M-2$ and $m>m'$, 
and ${\bf r}^{(-)}_{{y}_z}[1]=[\dots,{r}^{(m,m')}_{{y}_z}[1],\dots]^T$ for $m=0,1,\dots,M-2$ and $m'>m$.
 
Let us 
consider an FC collecting the correlation vectors 
\vspace{-1mm}
\begin{equation}
{\bf r}_{{y}_{z}}=[{{\bf r}^{(0)T}_{{y}_{z}}[0]},{{\bf r}^{(+)T}_{{y}_{z}}[0]},{\bf r}^{(-)T}_{{y}_{z}}[1]]^T 
\label{eq:ryz}
\vspace{-1mm}
\end{equation}
from the $(z+1)$-th group of sensors, for $z=0,1,\dots,Z-1$. 
We collect $\{{\bf r}^{(0)}_{{y}_z}[0]\}_{z=0}^{Z-1}$ in~\eqref{eq:ryz} into ${\bf r}^{(0)}_y[0]=[{\bf r}^{(0)T}_{{y}_0}[0],$ ${\bf r}^{(0)T}_{{y}_1}[0],\dots,{\bf r}^{(0)T}_{{y}_{Z-1}}[0]]^T$, similarly define ${{\bf r}^{(+)}_y[0]}$ and ${\bf r}^{(-)}_y[1]$,  
and use 
\eqref{eq:rypzall} to write
\begin{subequations}\label{eq:ryasRcrxsimpler} 
\setlength{\abovedisplayskip}{3pt}
\setlength{\belowdisplayskip}{3pt}
\begin{align}
\label{eq:ryasRcrxsimpler0}{{\bf r}^{(0)}_y[0]}&={\bf r}_{c}[0]{r}_x[0]={\bf 1}_{MZ}{r}_x[0], \\
\begin{bmatrix} {{\bf r}^{(+)}_y[0]} \\
{\bf r}^{(-)}_y[1]
\end{bmatrix}&=
\begin{bmatrix}{\bf R}_{c}[-1]\\
{\bf R}_{c}[1]
\end{bmatrix}{\bf r}_x[1]={\bf R}_c{\bf r}_x[1],\label{eq:ryasRcrxsimpler1}
\end{align}
\end{subequations}
with ${\bf r}_{c}[0]=[{\bf r}^T_{c_0}[0],{\bf r}^T_{c_1}[0],\dots,{\bf r}^T_{c_{Z-1}}[0]]^T$ and with
${\bf R}_{c}[1]$ and ${\bf R}_{c}[-1]$ similarly defined as ${\bf r}_{c}[0]$. 
The FC can then use least-squares (LS) to reconstruct ${\bf r}_x[1]$ from ${{\bf r}^{(+)}_y[0]}$ and ${\bf r}^{(-)}_y[1]$ if the $(M^2-M)Z\times (N-1)$ matrix ${\bf R}_{c}$ in~\eqref{eq:ryasRcrxsimpler1} has full column rank. Meanwhile, it is clear from~\eqref{eq:ryasRcrxsimpler0} that ${r}_x[0]$ can always be reconstructed using LS and its value is given by the average of the entries of ${{\bf r}^{(0)}_y[0]}$. Defining ${\bf F}_{2N-1}$ as the $(2N-1)\times (2N-1)$ discrete Fourier transform matrix, the FC can then compute ${\bf r}_x[-1]$ from ${\bf r}_x[1]$ using the Hermitian symmetry of ${r}_x[\tilde{n}]$ and the power spectrum from the reconstructed ${\bf r}_{{x}}=[{r}_{{x}}[0],{\bf r}^T_{{x}}[1],{\bf r}^T_{{x}}[-1]]^T$ 
as
\vspace{-1.5mm}
\begin{equation}
{\bf p}_x={\bf F}_{2N-1}{\bf r}_x.
\label{eq:px_correlogram}
\vspace{-1.5mm}
\end{equation}

Due to the finite sensing time, all the correlations must be approximated from a finite number of samples. Consider the unbiased estimate of
${r}^{(m,m')}_{{y}_z}[l]$ in~\eqref{eq:rypzmm'filter} given by
\vspace{-1.5mm}
\begin{equation}
\hat{r}^{(m,m')}_{{y}_z}[l]=\frac{1}{P(L-|l|)}\sum_{p=0}^{P-1}\sum_{l'=\text{max}(0,l)}^{L-1+\text{min}(0,l)}{y}^{(m)}_{z,p}[l']{{y}^{(m')*}_{z,p}[l'-l]}, 
\label{eq:ryzmm'hat}
\end{equation}
for $m,m'=0,1,\dots,M-1$. Using~\eqref{eq:ryzmm'hat},
we then 
form the estimate of ${\bf r}_{{y}_{z}}$ in~\eqref{eq:ryz} as
\begin{equation}
\hat{\bf r}_{{y}_{z}}=[\hat{\bf r}^{(0)T}_{{y}_z}[0],\hat{\bf r}^{(+)T}_{y_z}[0],\hat{\bf r}^{(-)T}_{y_z}[1]]^T.
\label{eq:ryzhat}
\end{equation}
Next, we use~\eqref{eq:ryzhat} to form $\hat{\bf r}^{(0)T}_y[0]=[\hat{\bf r}^{(0)T}_{{y}_0}[0],\hat{\bf r}^{(0)T}_{{y}_1}[0],\dots,$ $\hat{\bf r}^{(0)T}_{{y}_{Z-1}}[0]]^T$ and, similarly also, 
${\hat{\bf r}^{(+)}_y[0]}$ and ${\hat{\bf r}^{(-)}_y[1]}$, 
as the estimates of 
${{\bf r}^{(0)}_y[0]}$, ${{\bf r}^{(+)}_y[0]}$, and ${\bf r}^{(-)}_y[1]$ 
in~\eqref{eq:ryasRcrxsimpler}, respectively. In this case, we apply LS on $\hat{\bf r}^{(0)}_y[0]$, ${\hat{\bf r}^{(+)}_y[0]}$, and ${\hat{\bf r}^{(-)}_y[1]}$ 
instead of ${{\bf r}^{(0)}_y[0]}$, ${{\bf r}^{(+)}_y[0]}$, and ${\bf r}^{(-)}_y[1]$ in~\eqref{eq:ryasRcrxsimpler}.

\section{Condition for Least-Squares Reconstruction}\label{correlogram_reconstruct}

We now focus on the LS reconstruction of ${\bf r}_x[1]$ in~\eqref{eq:ryasRcrxsimpler1} and first evaluate the full column rank condition of ${\bf R}_{c}$ for $Z=1$. 
Note from~\eqref{eq:rxpz_one} and~\eqref{eq:ryasRcrxsimpler1} 
that ${\bf R}_{c}={\bf R}_{c_z}$ for $Z=1$. 
We can find from~\eqref{eq:rczmm},~\eqref{eq:rczmlagmin1},~\eqref{eq:rczmlag1}, and~\eqref{eq:Rclags} 
that each row of ${\bf R}_{c_z}$ in~\eqref{eq:rxpz_one} has only a single one in one entry and zeros elsewhere. 
The full column rank condition of ${\bf R}_{c_z}$ is then ensured if and only if each of its columns has {\it at least} a single one. 
We define $\Omega(\mathcal{M}_z)=\{(n^{(m)}_{z}-n^{(m')}_{z})\text{ mod }N|\forall n^{(m)}_{z},n^{(m')}_{z}\in\mathcal{M}_z\}$ and review the definition of a circular sparse ruler in~\cite{Romero}, 
which is 
also known as a cyclic difference set~\cite{Miller}. 
\newline
\hspace*{1mm}{\it Definition 1: A length-$(N-1)$ circular sparse ruler is 
a set $\mathcal{Q} \subset \{0,1,\dots,N-1\}$ such that 
$\Omega(\mathcal{Q})=\{0,1,\dots,N-1\}$.} 
\newline
\hspace*{1mm}A length-$(N-1)$ circular sparse ruler is a ruler that 
can measure all integers 
from $0$ to $N-1$ in a modular fashion despite missing some of its integer marks.
Let us consider the following theorem. 
\newline
\hspace*{1mm}{\it Theorem 1: ${\bf R}_{c_z}$ in~\eqref{eq:rxpz_one} has full column rank if and only if 
$\mathcal{M}_z$ in~\eqref{eq:mathcalMz} 
is 
a circular sparse ruler of length $N-1$.}
\newline
\hspace*{1mm}{\it Proof:}  
Assuming $n^{(m)}_z>n^{(m')}_z$, having a one in the $(n^{(m)}_z-n^{(m')}_z)$-th entry of ${\bf r}^{(m,m')}_{c_z}[-1]$ in~\eqref{eq:rczmlagmin1} implies having a one in the $(n^{(m')}_z-n^{(m)}_z+N)$-th entry of ${\bf r}^{(m',m)}_{c_z}[1]$ in~\eqref{eq:rczmlag1}. 
It is then clear from~\eqref{eq:Rclags} that if ${\bf R}_{c_z}$ in~\eqref{eq:rxpz_one} has at least a single one in the $((n^{(m)}_z-n^{(m')}_z) \text{ mod }N)$-th column, it 
also has at least a single one in the $((n^{(m')}_z-n^{(m)}_z) \text{ mod }N)$-th column. 
As we need to ensure each column of 
${\bf R}_{c_z}$ to have at least a single one, ensuring ${\bf R}_{c_z}$ in~\eqref{eq:rxpz_one} to have full column rank 
is identical to 
ensuring that $\Omega(\mathcal{M}_z)=\{0,1,\dots,N-1\}$. 
By considering Definition~1, 
the proof is concluded. $\square$

The correlation reconstruction with $Z=1$ is equivalent to the work in~\cite{TSP12}, 
though~\cite{TSP12} does not directly 
relate 
the rank condition of the system matrix to 
a circular sparse ruler. 
For $Z>1$, 
observe from~\eqref{eq:ryasRcrxsimpler1} that ${\bf R}_{c_z}$ in~\eqref{eq:rxpz_one} 
does not need to have full column rank and thus 
each of the sets $\{\mathcal{M}_z\}_{z=0}^{Z-1}$ does not have to be a circular sparse ruler. In fact, 
for $Z>1$, 
${\bf R}_{c}$ in~\eqref{eq:ryasRcrxsimpler1} has full column rank 
if 
\vspace{-2mm}
\begin{equation}
\bigcup_{z=0}^{Z-1}\Omega(\mathcal{M}_z)=\{0,1,\dots,N-1\}.
\label{eq:bigcupMz}
\vspace{-2mm}
\end{equation}
For completeness, let us consider the following definition.
\newline
\hspace*{1mm}{\it Definition 2: A length-$(N-1)$ incomplete circular sparse ruler 
is 
a set $\mathcal{Q}\subset \{0,1,\dots,N-1\}$ such that $\Omega(\mathcal{Q}) \subsetneq \{0,1,\dots,N-1\}$.}
\newline
\hspace*{1mm}Like a length-$(N-1)$ circular sparse ruler, 
a length-$(N-1)$ incomplete circular sparse ruler 
misses some of its integer marks, but 
it cannot measure all integers 
from $0$ to $N-1$ in a modular fashion.
Using Definition~2, we can view the design of $\{\mathcal{M}_z\}_{z=0}^{Z-1}$ satisfying~\eqref{eq:bigcupMz} as having $Z$ circular rulers, preferably sparse and incomplete to achieve 
a strong compression, such that, for each 
$n\in \{0,1,\dots,N-1\}$, 
at least one of these $Z$ incomplete circular sparse rulers 
can measure the distance 
$n$ in a modular fashion. 
Clearly, there is 
a tradeoff between $Z$ and $M$ under the constraint of~\eqref{eq:bigcupMz} since it is not possible to minimize 
both of them. 
Here, we aim at minimizing $Z$ given $M$, i.e., 
\vspace{-1mm}
\begin{equation}
\min_{Z,\{\mathcal{M}_z\}_{z=0}^{Z-1}} Z \text{ s.t. } \eqref{eq:bigcupMz}
\text{ and }|\mathcal{M}_z|=M, \:\forall z.
\label{eq:GivenMminZ}
\vspace{-1mm}
\end{equation}
Note that maximizing $|\Omega(\mathcal{M}_z)|$ in~\eqref{eq:bigcupMz} for all $z$ implies minimizing $Z$. 
As we have $((n^{(m')}_z-n^{(m)}_z) \text{ mod } N) \in \Omega(\mathcal{M}_z)$ and $((n^{(m)}_z-n^{(m')}_z) \text{ mod } N) \in \Omega(\mathcal{M}_z)$
whenever $n^{(m)}_z,n^{(m')}_z \in \mathcal{M}_z$,   
we can find that $|\Omega(\mathcal{M}_z)|\leq M(M-1)+1$. 
Disregarding the self-difference of the coset indices, we can find from~\eqref{eq:bigcupMz} that $Z$ is 
bounded as $Z\geq\left\lceil\frac{N-1}{M(M-1)}\right\rceil$. 
We thus aim to have $|\Omega(\mathcal{M}_z)|$ as close as possible to $M(M-1)+1$ for all $z$. 
\newline
\hspace*{1mm}{\it Definition~3: A length-$(N-1)$ circular Golomb ruler 
is 
a set $\mathcal{Q}\subset \{0,1,\dots,N-1\}$ such that, if 
$q,q',\tilde{q},\tilde{q}'\in\mathcal{Q}$ with $q\neq q'$, then $(q-q') \text{ mod N }=(\tilde{q}-\tilde{q}')\text{ mod N }$ implies $q=\tilde{q}$ and $q'=\tilde{q}'$~\cite{ColbournDinitz}. 
We say that 
$Z$ length-$(N-1)$ circular Golomb rulers $\{\mathcal{Q}_z\}_{z=0}^{Z-1}$ 
are non-overlapping if 
$\bigcap_{z=0}^{Z-1}\Omega(\mathcal{Q}_z)=\{0\}$.}
\newline
\hspace*{1mm}
Definition~3 implies that $|\Omega(\mathcal{M}_z)|=M(M-1)+1$ if and only if $\mathcal{M}_z$ is a circular Golomb ruler. Hence, a 
way to (approximately) minimize $Z$ given $M$ is to search for 
{\it non-overlapping} circular Golomb rulers that cover 
all or a certain number of the $N$ integer distances. 
In the latter case, we continue 
to search for another circular Golomb or an ordinary incomplete circular sparse ruler with $M$ marks that covers as many of the remaining uncovered integer distances as possible. This 
step is repeated until all the $N$ integer distances are covered. 
For $M=2$, the lower bound for $Z$ is 
$\left\lceil\frac{N-1}{2}\right\rceil$, which is 
reached for $N$ odd by having 
$\frac{N-1}{2}$ non-overlapping circular Golomb rulers with the $(z+1)$-th ruler $\mathcal{M}_{z}=\{0,z+1\}$.
This bound is also reached for $N$ even by having $\frac{N}{2}-1$ non-overlapping circular Golomb rulers with 
$\mathcal{M}_{z}=\{0,z+1\}$,
plus one incomplete circular sparse ruler $\mathcal{M}_{\frac{N}{2}-1}=\{0,\frac{N}{2}\}$ as the last ruler.
For $M=3$, the lower bound for $Z$ is 
$\left\lceil\frac{N-1}{6}\right\rceil$.
Table~\ref{tab:circularGolombruler} provides a list of non-overlapping circular Golomb rulers that achieve this bound for integers $\frac{N-1}{6}$ and $43\leq N\leq 115$. These non-overlapping circular Golomb rulers are found by first having $\mathcal{M}_z=\{0,z+1\}$, for all $z=0,1,\dots,Z-1$. The third (last) coset index in $\mathcal{M}_z$, for each $z$, is then determined by using a greedy search from the remaining coset indices that have not yet been used in $\{\mathcal{M}_z\}_{z=0}^{Z-1}$.
We are still investigating if there is {\it always} a sampling pattern that reaches the lower bound for $Z$, for any value of $N$ and $M$, and if there is a better algorithm to find this optimal sampling pattern. 
A set of non-overlapping {\it linear} (instead of circular) Golomb rulers called a {\it perfect difference basis system} is discussed in~\cite{Wild_Difference_Basis}. 
\section{Numerical Study}\label{numerical}

\begin{figure}[ht]
				\centering
        \includegraphics[width=0.48\textwidth]{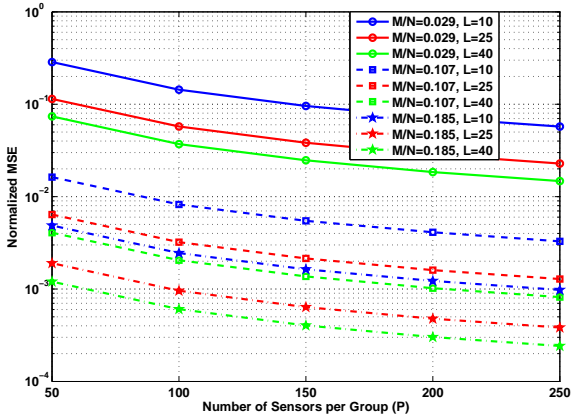}
        \vspace{-3mm}
        \caption{The NMSE between the compressively reconstructed power spectrum and the one reconstructed from the Nyquist rate samples.}
        \label{fig:NormalizedMSE_PML}
\end{figure}

We consider $N=103$ and six user signals whose frequency bands are given in Table~\ref{tab:experiment1} together with the power at each band normalized by frequency. These signals are generated by passing six sets of circular complex zero-mean Gaussian i.i.d. noise signals, with the variances set according to the desired user signal powers, 
into different digital filters having $N$ taps where the location of the unit-gain passband of the filter for each realization corresponds to the six different active bands. We assume that the 
signals are observed by unsynchronized sensors, which means that, at a certain point in time, all sensors generally observe different parts of the user signals. To simplify the experiment, we assume that, at time $t$, the $(zP+p)$-th sensor observes the part of the user signals that has previously been observed by the $(zP+p-1)$-th sensor at time $t-14T$, with $T$ the Nyquist sampling time. 
We start from $M=3$ active cosets per sensor and fix $Z$ at its lower bound, which is $Z=17$ for 
$M=3$, 
by having $17$ non-overlapping circular Golomb rulers (search for $N=103$ in Table~\ref{tab:circularGolombruler}).
We set the white noise power at each sensor to 
$16\text{ dBm}$ and vary $M,P$ and $L$ (see Fig.~\ref{fig:NormalizedMSE_PML}). 
Each user signal received by different sensors is assumed to pass through different wireless channels but 
the signal from a user received by all sensors 
is assumed to experience the same path loss and shadowing. Table~\ref{tab:experiment1} indicates the amount of path loss experienced between each user and all sensors, which is assumed to include shadowing. We simulate the Rayleigh fading by generating the channel frequency response based on 
a zero-mean complex Gaussian distribution with variance governed by the path loss in Table~\ref{tab:experiment1}. 
Each band is assumed to experience flat fading. 
We compute the normalized mean square error (NMSE) of the compressively estimated power spectrum with respect to 
the Nyquist-rate based estimate (obtained by activating all ${N}$ cosets in each sensor). 
Fig.~\ref{fig:NormalizedMSE_PML} shows how increasing the compression rate per sensor up to 
$M/N=0.19$ 
improves the estimation quality. 
The estimation quality is also improved by having either more sensors $PZ$ or more samples per coset $L$. 
\vspace{-4mm}
\vspace{-2.67mm}
\begin{table}[ht]
	\caption{
	The frequency bands occupied by the users, their power, and the experienced path loss} 
	\centering
	\vspace{-1mm}
		\begin{tabular}{| c | c | c |}
			\hline
    User band (rad/sample) & Power/freq. (per rad/sample) &Path loss\\ \hline
    $[\frac{-8\pi}{9},\frac{-7\pi}{9}]$ & $38$ dBm& $-18$ dB\\ \hline
    $[\frac{-6\pi}{9},\frac{-5\pi}{9}]$ & $40$ dBm& $-19$ dB\\ \hline
    $[\frac{\pi}{9},\frac{2\pi}{9}]$ & $34$ dBm& $-11$ dB\\ \hline
    $[\frac{3\pi}{9},\frac{4\pi}{9}]$ & $34$ dBm& $-17$ dB\\ \hline
    $[\frac{4\pi}{9},\frac{5\pi}{9}]$ & $32$ dBm& $-13$ dB\\ \hline
    $[\frac{6\pi}{9},\frac{7\pi}{9}]$ & $35$ dBm& $-19$ dB\\ \hline
    \end{tabular}
\label{tab:experiment1}
\end{table}

\vspace{-0mm}
\begin{table}[ht]
	\caption{List of non-overlapping circular Golomb rulers that cover all of the $N$ integer distances} 
	\centering
		\begin{tabular}{| c | l |}
			\hline
    $N$ & Non-overlapping circular Golomb rulers\\ \hline
    43 & $\mathcal{M}_0=\{0,1,17\}$, $\mathcal{M}_1=\{0,2,12\}$, $\mathcal{M}_2=\{0,3,24\}$,\\
    	 & $\mathcal{M}_3=\{0,4,13\}$, $\mathcal{M}_4=\{0,5,28\}$, $\mathcal{M}_5=\{0,6,14\}$,\\ 
    	 & $\mathcal{M}_6=\{0,7,32\}$ \\\hline
    49 & $\mathcal{M}_0=\{0,1,13\}$, $\mathcal{M}_1=\{0,2,20\}$, $\mathcal{M}_2=\{0,3,14\}$,\\
    	 & $\mathcal{M}_3=\{0,4,30\}$, $\mathcal{M}_4=\{0,5,15\}$, $\mathcal{M}_5=\{0,6,27\}$,\\ 
    	 & $\mathcal{M}_6=\{0,7,16\}$, $\mathcal{M}_7=\{0,8,25\}$\\\hline
    55 & $\mathcal{M}_0=\{0,1,15\}$, $\mathcal{M}_1=\{0,2,23\}$, $\mathcal{M}_2=\{0,3,16\}$,\\
    	 & $\mathcal{M}_3=\{0,4,33\}$, $\mathcal{M}_4=\{0,5,17\}$, $\mathcal{M}_5=\{0,6,31\}$,\\ 
    	 & $\mathcal{M}_6=\{0,7,18\}$, $\mathcal{M}_7=\{0,8,28\}$, $\mathcal{M}_8=\{0,9,19\}$\\\hline
    61 & $\mathcal{M}_0=\{0,1,16\}$, $\mathcal{M}_1=\{0,2,39\}$, $\mathcal{M}_2=\{0,3,17\}$,\\
    	 & $\mathcal{M}_3=\{0,4,38\}$, $\mathcal{M}_4=\{0,5,18\}$, $\mathcal{M}_5=\{0,6,35\}$,\\ 
    	 & $\mathcal{M}_6=\{0,7,19\}$, $\mathcal{M}_7=\{0,8,36\}$, $\mathcal{M}_8=\{0,9,20\}$,\\
    	 & $\mathcal{M}_9=\{0,10,31\}$\\ \hline
    67 & $\mathcal{M}_0=\{0,1,18\}$, $\mathcal{M}_1=\{0,2,42\}$, $\mathcal{M}_2=\{0,3,19\}$,\\
    	 & $\mathcal{M}_3=\{0,4,41\}$, $\mathcal{M}_4=\{0,5,20\}$, $\mathcal{M}_5=\{0,6,35\}$,\\ 
    	 & $\mathcal{M}_6=\{0,7,21\}$, $\mathcal{M}_7=\{0,8,36\}$, $\mathcal{M}_8=\{0,9,22\}$,\\
    	 & $\mathcal{M}_9=\{0,10,34\}$, $\mathcal{M}_{10}=\{0,11,23\}$\\ \hline
    73 & $\mathcal{M}_0=\{0,1,43\}$, $\mathcal{M}_1=\{0,2,20\}$, $\mathcal{M}_2=\{0,3,44\}$,\\
    	 & $\mathcal{M}_3=\{0,4,21\}$, $\mathcal{M}_4=\{0,5,40\}$, $\mathcal{M}_5=\{0,6,22\}$,\\ 
    	 & $\mathcal{M}_6=\{0,7,34\}$, $\mathcal{M}_7=\{0,8,23\}$, $\mathcal{M}_8=\{0,9,28\}$,\\
    	 & $\mathcal{M}_{9}=\{0,10,24\}$, $\mathcal{M}_{10}=\{0,11,37\}$, $\mathcal{M}_{11}=\{0,12,25\}$\\ \hline
    79 & $\mathcal{M}_0=\{0,1,28\}$, $\mathcal{M}_1=\{0,2,21\}$, $\mathcal{M}_2=\{0,3,48\}$,\\
    	 & $\mathcal{M}_3=\{0,4,22\}$, $\mathcal{M}_4=\{0,5,40\}$, $\mathcal{M}_5=\{0,6,23\}$,\\ 
    	 & $\mathcal{M}_6=\{0,7,49\}$, $\mathcal{M}_7=\{0,8,24\}$, $\mathcal{M}_8=\{0,9,50\}$,\\
    	 & $\mathcal{M}_9=\{0,10,25\}$, $\mathcal{M}_{10}=\{0,11,47\}$, $\mathcal{M}_{11}=\{0,12,26\}$,\\ 
    	 & $\mathcal{M}_{12}=\{0,13,33\}$\\\hline
		85 & $\mathcal{M}_0=\{0,1,32\}$, $\mathcal{M}_1=\{0,2,23\}$, $\mathcal{M}_2=\{0,3,50\}$,\\
    	 & $\mathcal{M}_3=\{0,4,24\}$, $\mathcal{M}_4=\{0,5,51\}$, $\mathcal{M}_5=\{0,6,25\}$,\\ 
    	 & $\mathcal{M}_6=\{0,7,48\}$, $\mathcal{M}_7=\{0,8,26\}$, $\mathcal{M}_8=\{0,9,49\}$,\\
    	 & $\mathcal{M}_9=\{0,10,27\}$, $\mathcal{M}_{10}=\{0,11,33\}$, $\mathcal{M}_{11}=\{0,12,28\}$,\\ 
    	 & $\mathcal{M}_{12}=\{0,13,43\}$, $\mathcal{M}_{13}=\{0,14,29\}$\\\hline
		91 & $\mathcal{M}_0=\{0,1,59\}$, $\mathcal{M}_1=\{0,2,24\}$, $\mathcal{M}_2=\{0,3,54\}$,\\
    	 & $\mathcal{M}_3=\{0,4,25\}$, $\mathcal{M}_4=\{0,5,41\}$, $\mathcal{M}_5=\{0,6,26\}$,\\ 
    	 & $\mathcal{M}_6=\{0,7,56\}$, $\mathcal{M}_7=\{0,8,27\}$, $\mathcal{M}_8=\{0,9,52\}$,\\
    	 & $\mathcal{M}_9=\{0,10,28\}$, $\mathcal{M}_{10}=\{0,11,57\}$, $\mathcal{M}_{11}=\{0,12,29\}$,\\ 
    	 & $\mathcal{M}_{12}=\{0,13,60\}$, $\mathcal{M}_{13}=\{0,14,30\}$, $\mathcal{M}_{14}=\{0,15,38\}$\\\hline
		97 & $\mathcal{M}_0=\{0,1,18\}$, $\mathcal{M}_1=\{0,2,21\}$, $\mathcal{M}_2=\{0,3,23\}$,\\
    	 & $\mathcal{M}_3=\{0,4,26\}$, $\mathcal{M}_4=\{0,5,29\}$, $\mathcal{M}_5=\{0,6,31\}$,\\ 
    	 & $\mathcal{M}_6=\{0,7,34\}$, $\mathcal{M}_7=\{0,8,36\}$, $\mathcal{M}_8=\{0,9,65\}$,\\
    	 & $\mathcal{M}_9=\{0,10,60\}$, $\mathcal{M}_{10}=\{0,11,46\}$, $\mathcal{M}_{11}=\{0,12,67\}$,\\ 
    	 & $\mathcal{M}_{12}=\{0,13,52\}$, $\mathcal{M}_{13}=\{0,14,54\}$, $\mathcal{M}_{14}=\{0,15,59\}$,\\
    	 & $\mathcal{M}_{15}=\{0,16,64\}$\\\hline
   103 & $\mathcal{M}_0=\{0,1,27\}$, $\mathcal{M}_1=\{0,2,39\}$, $\mathcal{M}_2=\{0,3,28\}$,\\
    	 & $\mathcal{M}_3=\{0,4,62\}$, $\mathcal{M}_4=\{0,5,29\}$, $\mathcal{M}_5=\{0,6,65\}$,\\ 
    	 & $\mathcal{M}_6=\{0,7,30\}$, $\mathcal{M}_7=\{0,8,57\}$, $\mathcal{M}_8=\{0,9,31\}$,\\
    	 & $\mathcal{M}_9=\{0,10,50\}$, $\mathcal{M}_{10}=\{0,11,32\}$, $\mathcal{M}_{11}=\{0,12,60\}$,\\ 
    	 & $\mathcal{M}_{12}=\{0,13,33\}$, $\mathcal{M}_{13}=\{0,14,56\}$, $\mathcal{M}_{14}=\{0,15,34\}$,\\
    	 & $\mathcal{M}_{15}=\{0,16,52\}$, $\mathcal{M}_{16}=\{0,17,35\}$\\\hline   
   109 & $\mathcal{M}_0=\{0,1,28\}$, $\mathcal{M}_1=\{0,2,42\}$, $\mathcal{M}_2=\{0,3,29\}$,\\
    	 & $\mathcal{M}_3=\{0,4,61\}$, $\mathcal{M}_4=\{0,5,30\}$, $\mathcal{M}_5=\{0,6,70\}$,\\ 
    	 & $\mathcal{M}_6=\{0,7,31\}$, $\mathcal{M}_7=\{0,8,68\}$, $\mathcal{M}_8=\{0,9,32\}$,\\
    	 & $\mathcal{M}_9=\{0,10,66\}$, $\mathcal{M}_{10}=\{0,11,33\}$, $\mathcal{M}_{11}=\{0,12,71\}$,\\ 
    	 & $\mathcal{M}_{12}=\{0,13,34\}$, $\mathcal{M}_{13}=\{0,14,58\}$, $\mathcal{M}_{14}=\{0,15,35\}$,\\
    	 & $\mathcal{M}_{15}=\{0,16,63\}$, $\mathcal{M}_{16}=\{0,17,36\}$, $\mathcal{M}_{17}=\{0,18,55\}$\\\hline   
   115 & $\mathcal{M}_0=\{0,1,30\}$, $\mathcal{M}_1=\{0,2,45\}$, $\mathcal{M}_2=\{0,3,31\}$,\\
    	 & $\mathcal{M}_3=\{0,4,64\}$, $\mathcal{M}_4=\{0,5,32\}$, $\mathcal{M}_5=\{0,6,73\}$,\\ 
    	 & $\mathcal{M}_6=\{0,7,33\}$, $\mathcal{M}_7=\{0,8,71\}$, $\mathcal{M}_8=\{0,9,34\}$,\\
    	 & $\mathcal{M}_9=\{0,10,69\}$, $\mathcal{M}_{10}=\{0,11,35\}$, $\mathcal{M}_{11}=\{0,12,74\}$,\\ 
    	 & $\mathcal{M}_{12}=\{0,13,36\}$, $\mathcal{M}_{13}=\{0,14,61\}$, $\mathcal{M}_{14}=\{0,15,37\}$,\\
    	 & $\mathcal{M}_{15}=\{0,16,66\}$, $\mathcal{M}_{16}=\{0,17,38\}$, $\mathcal{M}_{17}=\{0,18,58\}$,\\
    	 & $\mathcal{M}_{18}=\{0,19,39\}$\\\hline   
     \end{tabular}
\label{tab:circularGolombruler}
\end{table}
\end{document}